%

\catcode`\@=11

%
%
\def\bibn@me{R\'ef\'erences}
\def\bibliographym@rk{\centerline{{\sc\bibn@me}}
	\sectionmark\section{\ignorespaces}{\unskip\bibn@me}
	\bigbreak\bgroup
	\ifx\ninepoint\undefined\relax\else\ninepoint\fi}
%
%
%
\let\refsp@ce=\
\let\bibleftm@rk=[
\let\bibrightm@rk=]
%
%
%
\def\numero{n\raise.82ex\hbox{$\fam0\scriptscriptstyle
o$}~\ignorespaces}
%
%
\newcount\equationc@unt
\newcount\bibc@unt
\newif\ifref@changes\ref@changesfalse
\newif\ifpageref@changes\ref@changesfalse
\newif\ifbib@changes\bib@changesfalse
\newif\ifref@undefined\ref@undefinedfalse
\newif\ifpageref@undefined\ref@undefinedfalse
\newif\ifbib@undefined\bib@undefinedfalse
\newwrite\@auxout
%
%
%
%
%
%
%
%
\def\re@dreferences#1#2{{%
	\re@dreferenceslist{#1}#2,\undefined\@@}}
\def\re@dreferenceslist#1#2,#3\@@{\def\next{#2}%
	\expandafter\ifx\csname#1@@\meaning\next\endcsname\relax
	??\immediate\write16
	{Warning, #1-reference "\next" on page \the\pageno\space
	is undefined.}%
	\global\csname#1@undefinedtrue\endcsname
	\else\csname#1@@\meaning\next\endcsname\fi
	\ifx#3\undefined\relax
	\else,\refsp@ce\re@dreferenceslist{#1}#3\@@\fi}
%
%
%
\def\newlabel#1#2{{\def\next{#1}\newl@bel#2}}
\def\newl@bel#1#2{%
	\expandafter\xdef\csname ref@@\meaning\next\endcsname{#1}%
	\expandafter\xdef\csname pageref@@\meaning\next\endcsname{#2}}
\def\label#1{{%
	\toks0={#1}\message{ref(\lastref) \the\toks0,}%
	\ignorespaces\immediate\write\@auxout%
	{\noexpand\newlabel{\the\toks0}{{\lastref}{\the\pageno}}}%
	\def\next{#1}%
	\expandafter\ifx\csname ref@@\meaning\next\endcsname\lastref%
	\else\global\ref@changestrue\fi%
	\newlabel{#1}{{\lastref}{\the\pageno}}}}
\def\ref#1{\re@dreferences{ref}{#1}}
\def\pageref#1{\re@dreferences{pageref}{#1}}
%
%
\def\bibcite#1#2{{\def\next{#1}%
	\expandafter\xdef\csname bib@@\meaning\next\endcsname{#2}}}
\def\cite#1{\bibleftm@rk\re@dreferences{bib}{#1}\bibrightm@rk}
%
%
\def\beginthebibliography#1{\bibliographym@rk
	\setbox0\hbox{\bibleftm@rk#1\bibrightm@rk\enspace}
	\parindent=\wd0
	\global\bibc@unt=0
	\def\bibitem##1{\global\advance\bibc@unt by 1
		\edef\lastref{\number\bibc@unt}
		{\toks0={##1}
		\message{bib[\lastref] \the\toks0,}%
		\immediate\write\@auxout
		{\noexpand\bibcite{\the\toks0}{\lastref}}}
		\def\next{##1}%
		\expandafter\ifx
		\csname bib@@\meaning\next\endcsname\lastref
		\else\global\bib@changestrue\fi%
		\bibcite{##1}{\lastref}
		\medbreak
		\item{\hfill\bibleftm@rk\lastref\bibrightm@rk}%
		}
	}
\def\endthebibliography{\egroup\par}
%
%
    \outer\def\bye{\@closeaux
    	\par\vfill\supereject\end}
%
\def\@closeaux{\closeout\@auxout
	\ifref@changes\immediate\write16
	{Warning, changes in references.}\fi
	\ifpageref@changes\immediate\write16
	{Warning, changes in page references.}\fi
	\ifbib@changes\immediate\write16
	{Warning, changes in bibliography.}\fi
	\ifref@undefined\immediate\write16
	{Warning, references undefined.}\fi
	\ifpageref@undefined\immediate\write16
	{Warning, page references undefined.}\fi
	\ifbib@undefined\immediate\write16
	{Warning, citations undefined.}\fi}
%
%
\immediate\openin\@auxout=\jobname.aux
\ifeof\@auxout \immediate\write16
     {Creating file \jobname.aux}
\immediate\closein\@auxout
\immediate\openout\@auxout=\jobname.aux
\immediate\write\@auxout {\relax}%
\immediate\closeout\@auxout
\else\immediate\closein\@auxout\fi
%
%
\input\jobname.aux \par
\immediate\openout\@auxout=\jobname.aux
%
%

\def\bibn@me{R\'ef\'erences bibliographiques}
\catcode`@=11
\def\bibliographym@rk{\bgroup}
%
%
\outer\def\bye{ 	\par\vfill\supereject\end}

\magnification=1200

\font\tenbfit=cmbxti10 
\font\sevenbfit=cmbxti10 at 7pt 
\font\sixbfit=cmbxti5 at 6pt 

\newfam\mathboldit 

\textfont\mathboldit=\tenbfit
  \scriptfont\mathboldit=\sevenbfit
   \scriptscriptfont\mathboldit=\sixbfit

\def\bfit           
{\tenbfit           
   \fam\mathboldit 
}

\def\Z{{\bf Z}}     
\def\R{{\bf R}}   \def\cK{{\bf B}}

\def\cK{{\cal {K}}}

\def\cE{{\cal {E}}}

\def\Bad{{\bfit Bad}}

\def\house#1{\setbox1=\hbox{$\,#1\,$}%
\dimen1=\ht1 \advance\dimen1 by 2pt \dimen2=\dp1 \advance\dimen2 by 2pt
\setbox1=\hbox{\vrule height\dimen1 depth\dimen2\box1\vrule}%
\setbox1=\vbox{\hrule\box1}%
\advance\dimen1 by .4pt \ht1=\dimen1
\advance\dimen2 by .4pt \dp1=\dimen2 \box1\relax}

  \def\eps{{\varepsilon}}

\def\build#1_#2^#3{\mathrel{\mathop{\kern 0pt#1}\limits_{#2}^{#3}}}

\def\date {le\ {\the\day}\ \ifcase\month\or janvier
\or fevrier\or mars\or avril\or mai\or juin\or juillet\or
ao\^ut\or septembre\or octobre\or novembre
\or d\'ecembre\fi\ {\oldstyle\the\year}}

\font\fivegoth=eufm5 \font\sevengoth=eufm7 \font\tengoth=eufm10

\newfam\gothfam \scriptscriptfont\gothfam=\fivegoth
\textfont\gothfam=\tengoth \scriptfont\gothfam=\sevengoth

\def\pro{\noindent {\it Proof. }}

\def\smallsquare{\vbox{\hrule\hbox{\vrule height 1 ex\kern 1 ex\vrule}\hrule}}
\def\cqfd{\hfill \smallsquare\vskip 3mm}



\centerline{}

\vskip 4mm

\centerline{
\bf Badly approximable numbers and Littlewood-type problems}

\vskip 8mm
\centerline{Yann B{\sevenrm UGEAUD}
\ \& \ Nikolay M{\sevenrm OSHCHEVITIN}
\footnote{}{\rm 
2000 {\it Mathematics Subject Classification : 11J13; 11J25, 11K60} .}
}

{\narrower\narrower
\vskip 12mm

\proclaim Abstract. {
We establish that the set of pairs $(\alpha, \beta)$ of real numbers 
such that
$$
\liminf_{q \to + \infty} \, q \cdot 
(\log q)^2 \cdot \Vert q \alpha \Vert \cdot \Vert q \beta \Vert > 0, 
$$
where $\Vert \cdot \Vert$ denotes the distance to the 
nearest integer,
has full Hausdorff dimension in $\R^2$.
Our proof rests on a method introduced by Peres and Schlag,
that we further apply to various Littlewood-type problems.
}

}

\vskip 15mm

\centerline{\bf 1. Introduction}

\vskip 6mm

A famous open problem in simultaneous 
Diophantine approximation, called 
the Littlewood conjecture \cite{Lit68}, claims that, 
for any given pair $(\alpha, \beta)$ of real numbers, 
we have
$$
\inf_{q \ge 1} \, q \cdot \Vert q \alpha \Vert 
\cdot \Vert q \beta \Vert = 0,
\eqno (1.1)
$$
where $\Vert \cdot \Vert$ denotes the distance to the 
nearest integer.
Throughout the present paper, we denote by $\Bad$ the set of badly
approximable numbers, that is,
$$
\Bad = \{ \alpha \in \R : \inf_{q \ge 1} \,
q \cdot \Vert q \alpha \Vert > 0\},
$$
and we recall that $\Bad$ has Lebesgue measure zero
and full Hausdorff dimension \cite{Ja28}.
Consequently, (1.1) holds for almost every
pair $(\alpha, \beta)$ of real numbers.
Recently, this result was considerably improved by
Einsiedler, Katok and Lindenstrauss \cite{EKL},
who established that the set of pairs $(\alpha, \beta)$
for which (1.1) do not hold has Hausdorff dimension zero;
see also \cite{PoVe} for a weaker statement, 
and Section 10.1 of
\cite{BuLiv} for a survey of related results.

Another metrical statement connected
to the Littlewood conjecture was
established by Gallagher \cite{Gal62}
in 1962 and can be formulated as follows 
(see e.g. \cite{BeVe07}).

\proclaim Theorem G.
Let $n$ be a positive integer.
Let $\Psi : \R_{>0} \to \R_{>0}$ be a non-increasing
function.
The set of points $(x_1, \ldots , x_n)$ in $\R^n$
such that there are infinitely many positive integers $q$
satisfying
$$
\prod_{i=1}^n \, \Vert q x_i \Vert < \Psi (q)
$$
has full Lebesgue measure if the sum
$$
\sum_{h \ge 1} \, \Psi(h)^n \, (\log h)^{n-1}
$$
diverges, and has zero Lebesgue measure otherwise.

In particular, it follows from Gallagher's theorem that
$$
\liminf_{q \to + \infty} \, q \cdot 
(\log q)^2 \cdot \Vert q \alpha \Vert \cdot \Vert q \beta \Vert = 0
\eqno (1.2)
$$
for almost every pair $(\alpha, \beta)$ of real numbers.
The main purposes of the present note
are to establish the existence of
exceptional pairs $(\alpha, \beta)$ which do not satisfy (1.2)
--- a result first proved in \cite{Mo5} ---,
and to prove that the set of these pairs has
full Hausdorff dimension in $\R^2$.
We further consider various questions closely
related to the Littlewood conjecture.

Our main results are stated in Section 2 
and proved in Sections 4 and 5,
with the help of auxiliary lemmas gathered in Section 3.
Several additional results are given in Section~6.

Throughout this paper, $\lfloor x \rfloor$ and $\lceil x \rceil$
denote the greatest integer less than or equal to $x$
and the smallest integer greater than or equal to $x$, respectively.

\vskip 8mm

\centerline{\bf 2. Main results}

\vskip 6mm

Our first result shows that there are many pairs $(\alpha, \beta)$
of real numbers that are not well multiplicatively approximable.

\proclaim Theorem 1.
For every real number $\alpha$ in $\Bad$, the set of real numbers $\beta$
such that
$$
\liminf_{q \to + \infty} \, q \cdot 
(\log q)^2 \cdot \Vert q \alpha \Vert 
\cdot \Vert q \beta \Vert > 0  \eqno (2.1)
$$
has full Hausdorff dimension.

The proof of Theorem 1 uses a method 
introduced by Peres and Schlag \cite{PeSc},
which was subsequently applied in \cite{Mo1,Mo2,Mo3,Mo4,Mo5}.

Since the set $\Bad$ has full Hausdorff dimension,
the next result follows from Theorem~1 by an 
immediate application of
Corollary 7.12 from \cite{Fal90}.

\proclaim Theorem 2.
The set of pairs $(\alpha, \beta)$ of real numbers satisfying
$$
\liminf_{q \to + \infty} \, q \cdot 
(\log q)^2 \cdot \Vert q \alpha \Vert \cdot \Vert q \beta \Vert > 0  
$$
has full Hausdorff dimension in $\R^2$.

Theorem 1 can be viewed as a complement to the following
result of Pollington and Velani \cite{PoVe}.

\proclaim Theorem PV.
For every real number $\alpha$ in $\Bad$, there exists a
subset $G(\alpha)$ of $\Bad$ with full Hausdorff dimension
such that, for any $\beta$ in $G(\alpha)$, there
exist arbitrarily large integers $q$ satisfying
$$
q \cdot (\log q) \cdot \Vert q \alpha \Vert \cdot \Vert q \beta \Vert \le 1.
$$

In \cite{AdBu06}, the authors constructed explicitly for every $\alpha$
in $\Bad$ uncountably many $\beta$ in $\Bad$ such that the
pair $(\alpha, \beta)$ satisfies (1.1), and even a strong form of
this inequality. It would be very interesting to construct
explicit examples of pairs of real numbers that satisfy (2.1).

A modification of an auxiliary lemma
yields a slight improvement on Theorem~1.

\proclaim Theorem 3.
Let $a$ be a real number with $0 < a < 1$.
For every real number $\alpha$ in $\Bad$, the set of real numbers $\beta$
such that
$$
\liminf_{q \to + \infty} \, q \cdot 
(\log q)^{2-a} \cdot (\log 1/\Vert q \alpha \Vert)^a
\cdot \Vert q \alpha \Vert \cdot \Vert q \beta \Vert > 0  
$$
has full Hausdorff dimension.

Theorem 3 is stronger than Theorem 1 since,
for every $\alpha$ in $\Bad$, there exists
a positive real number $\delta$ such that
$\log (1/\Vert q \alpha \Vert) \le \delta \log q$ holds
for every integer $q \ge 2$.

Cassels and Swinnerton-Dyer \cite{CaSw} 
proved that (1.1)
is equivalent to the equality 
$$
\inf_{(x, y)\in \Z \times \Z \setminus\{(0,0)\}} \, 
\max\{\vert x \vert,1\}\cdot\max\{\vert y \vert,1\}\cdot \Vert x \alpha 
+ y \beta\Vert = 0,
$$ 
and used it 
to show that (1.1) holds if $\alpha$ and $\beta$ belong to the same
cubic number field (see also \cite{Peck}).
In this context, we have the following metrical result,
extracted from page 455 of \cite{BeKlMa}.
For integers $q_1, \ldots , q_n$, set
$$
\Pi(q_1, \ldots , q_n) = \prod_{i=1}^n \, \max\{1, \vert q_i \vert\}.
$$

\proclaim Theorem BKM.
Let $n$ be a positive integer.
Let $\Psi : \R_{>0} \to \R_{>0}$ be a non-increasing
function.
The set of points $(x_1, \ldots , x_n)$ in $\R^n$
such that there are infinitely many integers $q_1, \ldots, q_n$
satisfying
$$
|| q_1 x_1 + \ldots + q_n x_n || < \Psi \bigl( \Pi(q_1, \ldots , q_n) \bigr)
\eqno (2.2)
$$
has full Lebesgue measure if the sum
$$
\sum_{h \ge 1} \, \Psi(h) \, (\log h)^{n-1}   \eqno (2.3)
$$
diverges, and has zero Lebesgue measure otherwise.

For $n \ge 2$, there is no known example 
of points $(x_1, \ldots , x_n)$ in $\R^n$
and of a function $\Psi$ as in Theorem BKM such that the sum 
(2.3) diverges and (2.2) has only finitely many solutions.
The Peres--Schlag method allows us
to show that such examples do exist.

\proclaim Theorem 4.
The set of pairs $(\alpha, \beta)$ of real numbers satisfying
$$
\liminf_{x, y \ge 0} \, \max\{2, |xy|\} \cdot 
\Vert x \alpha + y \beta \Vert \cdot 
(\log  \max\{2, |xy|\})^2 > 0  
$$
has full Hausdorff dimension in $\R^2$.

The proof of Theorem 4 is briefly outlined in Section 5.
Note that Theorem 4 (resp. Theorem 1) 
does not follow from Theorem 1 (resp. Theorem 4)
by some transference principle.

In analogy with the Littlewood conjecture,
de Mathan and Teuli\'e \cite{BdMTe} 
proposed recently a `mixed Littlewood
conjecture'.
For any prime number $p$, the usual
$p$-adic absolute value $| \cdot |_p$ is
normalized in such a way that $|p|_p = p^{-1}$.

\proclaim De Mathan--Teuli\'e conjecture.
For every real number $\alpha$ and every prime
number $p$, we have
$$
\inf_{q \ge 1} \, q \cdot \Vert q \alpha \Vert \cdot
\vert q \vert_p = 0. 
$$

Despite several recent results \cite{EiKl07,BDM}, 
this conjecture is still unsolved.
The following metrical statement,
established in \cite{BuHaVe}, should be
compared with Theorem G.

\proclaim Theorem BHV.
Let $k$ be a positive integer.
Let $p_1, \ldots , p_k$ be distinct prime numbers.
Let $\Psi : \R_{>0} \to \R_{>0}$ be a non-increasing
function.
The set of real numbers $\alpha$
such that there are infinitely many positive integers $q$
satisfying
$$
\Vert q \alpha \Vert \cdot |q|_{p_1} \cdots 
|q|_{p_k} < \Psi (q)
$$
has full Lebesgue measure if the sum
$$
\sum_{h \ge 1} \, \Psi(h) \, (\log h)^k
$$
diverges, and has zero Lebesgue measure otherwise.

As an immediate consequence of Theorem BHV, we get that,
for every prime number $p$, almost every real number $\alpha$
satisfies
$$
\inf_{q \ge 2} \, q \cdot 
(\log q)^2 \cdot (\log \log q)
\cdot \Vert q \alpha \Vert \cdot \vert q \vert_p  =  0  \eqno (2.4)
$$
The method of proof of Theorem 1 allows us to confirm the
existence of real numbers for which (2.4) does not hold.

\proclaim Theorem 5.
Let $a$ be a real number with $0 \le a < 1$.
For every prime number $p$, the set of real numbers $\alpha$
such that
$$
\liminf_{q \to + \infty} \, q \cdot 
(\log q)^{2-a} \cdot \Vert q \alpha \Vert \cdot \vert q \vert_p 
\cdot (\log 2/ \vert q \vert_p)^a > 0  
$$
has full Hausdorff dimension.

We display an immediate consequence of Theorem 5.

\proclaim Corollary 1.
For every prime number $p$, the set of real numbers $\alpha$
such that
$$
\liminf_{q \to + \infty} \, q \cdot 
(\log q)^2 \cdot \Vert q \alpha \Vert \cdot \vert q \vert_p  > 0  
$$
has full Hausdorff dimension.

In the present note, we have restricted our attention to 
$2$-dimensional questions. However, our method can be successfully
applied to prove that, given an integer $n \ge 2$,
there are real numbers $\alpha_1, \ldots , \alpha_n$
such that
$$
\liminf_{q \to + \infty} \, q \cdot 
(\log q)^n \cdot \Vert q \alpha_1 \Vert \cdots \Vert q \alpha_n \Vert > 0,
$$
as well as real numbers $\beta_1, \ldots , \beta_n$
such that
$$
\liminf_{x_1, \ldots , x_n \ge 0} \, \max\{2, |x_1 \ldots x_n|\} \cdot 
\Vert x_1 \beta_1 + \ldots   + x_n \beta_n \Vert \cdot 
(\log  \max\{2, |x_1 \ldots x_n|\})^n > 0  
$$
This will be the subject of subsequent work by E. Ivanova.

\vskip 6mm

\centerline{\bf 3. Auxiliary results}

\vskip 4mm

The original method of Peres and Schlag is a construction of
nested intervals.
A useful tool for estimating from below the Hausdorff
measure of a Cantor set is the mass distribution principle, 
which we recall now.

We consider a set $\cK$ included in a bounded interval
$E$, and defined as follows.
Set $\cE_0 = E$ and assume that, for any positive integer $k$,
there exists a finite family $\cE_k$ of 
disjoint compact intervals in $E$ such that any interval $U$ belonging to
$\cE_k$ is contained in exactly
one of the intervals of $\cE_{k-1}$ and contains at least two intervals 
belonging to $\cE_{k+1}$. Suppose also that the
maximum of the lengths of the intervals in $\cE_k$ tends to 0
when $k$ tends to infinity. For $k \ge 0$, denote by $E_k$ 
the union of the intervals belonging to the family
$\cE_k$, and set
$$
\cK := \bigcap_{k=1}^{+ \infty} \, {E_k}.
$$

\proclaim Lemma 1. 
Keep the same notation as above.
Assume further that there exists 
a positive integer $k_0$ such that, for any
$k \ge k_0$, each interval of $E_{k-1}$ contains at least $m_k \ge 2$
intervals of $E_k$, these being separated by at least $\eps_k$, where
$0 < \eps_{k+1} < \eps_k$. We then have
$$
\dim \cK \ge \liminf_{k \to + \infty} \, 
{\log (m_1 \ldots m_{k-1}) \over - \log (m_k \eps_k)}. 
$$

\pro This is Example 4.6 in \cite{Fal90},
see also Proposition 5.2 in \cite{BuLiv}. \cqfd

\proclaim Lemma 2.
Let $\alpha$ be in $\Bad$.
There exists a positive constant $C(\alpha)$ such that, for every
integer $q \ge 2$, we have
$$
\sum_{x=q}^{q^3} \, {1 \over ||\alpha x||\, x\log_2^2 x} \le C(\alpha).
$$

\pro
This is a straightforward consequence of Example 3.2 on page 124
of \cite{KuNi}, where it is established that there exists 
a positive constant $C_1(\alpha)$ such that
$$
\sum_{x=1}^{m} \, {1 \over ||\alpha x||\, x} \le C_1(\alpha) (\log m)^2,
$$
for all positive integers $m$. \cqfd

Theorem 3 depends on the following refinement of Lemma 2.

\proclaim Lemma 3.
Let $\alpha$ be in $\Bad$.
Let $a$ be a real number with $0 < a < 1$.
There exists a positive constant $C(\alpha)$ such that, for every
integer $q \ge 2$, we have
$$
\sum_{x=q}^{q^3} \, {1 \over ||\alpha x||\, x \,
(\log 1/\Vert x \alpha \Vert)^a \cdot (\log x)^{2-a}} \le C(\alpha).
$$

\pro
Let $(p_j/q_j)_{j \ge 0}$ denote the sequence of convergents 
to $\alpha$. Let $m$ (resp. $n$) be the largest 
(resp. the smallest) integer $j$ such that $q_j \le q$
(resp. $q_j \ge q^3$).
As the sequence $(q_j)_{j \ge 0}$
grows exponentially fast, we have 
$$
\log q \ll n \ll m \ll \log q,
$$
where, as throughout this proof, the numerical constants
implied by $\ll$ depend only on $\alpha$.

Let $j$ be an integer satisfying $m \le j < n$
and consider 
$$
S_j := \sum_{x=q_j}^{q_{j+1}} \, {1 \over ||\alpha x||\, x \,
(\log 1/\Vert x \alpha \Vert)^a}.
$$
A classical result asserts that the points 
$\{ \alpha x \}$, $x = 1, \ldots , q_{j+1}$, are very well distributed
in $(0, 1)$. Consequently,
$$
\eqalign{
S_j \ll {1 \over q_j} \, 
 \sum_{x=1}^{q_{j+1}} \, {1 \over ||\alpha x||\, 
(\log 1/\Vert x \alpha \Vert)^a} 
& \ll {1 \over q_j} \, 
 \sum_{x=1}^{q_{j+1}/2} \, {(q_{j+1} /x) \over 
(\log (q_{j+1}/x))^a} \cr
& \ll {q_{j+1} \over q_j} \, \int_2^{q_{j+1}} \,
{du \over u (\log u)^{1 - a}} 
\ll (\log q_j)^{1 - a}, \cr}
$$
since $q_{j+1} / q_j$ is bounded from above by an absolute
constant depending only on $\alpha$. Now,
$$
\sum_{j=m}^n \, S_j \ll (\log q)^{2 - a},
$$
which proves the lemma. \cqfd

The key tool for the proof of Theorem 5 is Lemma 4 below.

\proclaim Lemma 4.
Let $p$ be a prime number.
Let $a$ be a real number with $0 \le a < 1$.
There exists a positive constant $C(a,p)$ such that, for every
integer $q \ge 2$, we have
$$
\sum_{x=q}^{q^3} \, {1 \over x \cdot |x|_p \, (\log (2/|x|_p))^a 
\cdot (\log x)^{2-a}} \le C(a,p).
$$

\pro
Observe that
$$
\sum_{x=q}^{q^3} \, {1 \over x \cdot |x|_p \cdot (\log (2/|x|_p))^a}
\ll \sum_{j=0}^{3 \log q} \,
\sum_{x= \lceil q/p^j \rceil}^{\lfloor q^3/p^j \rfloor} 
\, {1 \over x (j+1)^a},
$$
where the second summation is taken over the integers $x$
that are not divisible by $p$.
Consequently,
$$
\sum_{x=q}^{q^3} \, {1 \over x \cdot |x|_p \cdot (\log (2/|x|_p))^a}
\ll \sum_{j=0}^{3 \log q} \, {\log q \over (j+1)^a}
\ll (\log q)^{2-a},
$$
and the lemma is proved. \cqfd

\vskip 6mm

\centerline{\bf 4. Proof of Theorem 1}

\vskip 4mm

Let $\alpha$ be in $\Bad$ and $\delta$ 
be a positive real number satisfying
$$
q \cdot || q \alpha || \ge \delta, 
\quad \hbox{for every $q \ge 1$}. \eqno (4.1)
$$
Let $\eps$ be such that
$$
0 < \eps < \bigl(2^{10} C(\alpha) \bigr)^{-1}, \eqno (4.2)
$$
where $C(\alpha)$ is given by Lemma 2.

We follow a method introduced
by Peres and Schlag \cite{PeSc}.
First, we construct `dangerous' sets of real numbers. 
These sets depend
on $\alpha$, but, to simplify the notation, we choose not to
indicate this dependence.
For integers $x$ and $y$ with $x \ge 2$ and $0\le y\le x$, define
$$
E (x,y) = \biggl[ {y \over x} -{\varepsilon \over ||\alpha x || x^{2}\log_2^2 x}, 
{y \over x} 
+{\varepsilon \over ||\alpha x || x^{2}\log_2^2 x}
\biggr]  \eqno (4.3)
$$
and
$$
E (x) =\bigcup_{y =0}^{x} \, 
\bigl( \, E (x,y) \cap [0,1] \, \bigr). \eqno (4.4)
$$
Set also
$$
l_0 = 0, \quad
 l_x = \lfloor\log_2 (||\alpha x || x^{2}\log_2^2 x /(2\varepsilon) )  \rfloor ,
\quad \hbox{for $x \in \Z_{\ge 1}$}.  \eqno (4.5)
$$
Each interval from the union $E (x)$ defined in (4.4) 
can be covered by an open dyadic interval of the form
$$
\left( {b \over 2^{l_x }}, {b+z \over 2^{l_x }}\right),
\quad z = 1,2, \quad
b  \in \Z_{\ge 0}.
$$
Let $A (x)$ be the smallest union of all such dyadic intervals which 
covers the whole set $E (x)$ and put
$$
A^c (x) = [0,1] \setminus A (x).
$$ 
Observe that $A^c (x)$ is a union of closed intervals of the form
$$
\left[ {a \over 2^{l_x}}, {a+1 \over 2^{l_x}}\right] ,
\,\,\, a \in \Z_{\ge 0}. 
$$

Let $q_0$ be an integer such that
$$
q_0 \ge (100 \eps)^3  \quad 
\hbox{and} \quad
|| q_0 \alpha || \ge 1/4.  \eqno (4.6)
$$
For $q \ge q_0 $, define
$$
B_q =\bigcap_{x=q_0}^q A^c_{} (x).
$$ 
The sets $B_q$, $q \ge q_0$, are closed and nested.
Our aim is to show inductively that they are non-empty.
Set $L_0 = l_0$ and 
$$
q_k := q_0^{3^k}, \quad
L_k = \lfloor\log_2 ( q_k^2 \log_2^2 q_k / (4\varepsilon) )  \rfloor,
\quad k \ge 1.   \eqno (4.7)
$$
Observe that $l_x \le L_k$ when $x \le q_k$.

For every integer $k \ge 0$ we construct inductively subsets $C_{q_k}$
and $D_{q_k}$ of $B_{q_k}$ with the following property $(P_k)$:
\medskip
{\it The set
$C_{q_k}$ is the union of $2^{-5k-3 + L_k}$
intervals of length $2^{-L_k}$, separated by at least
$2^{-L_k}$, and such that at least $2^{-5k-5 + L_k}$
among them include at least $2^{L_{k+1}-L_k -3}$
intervals composing $B_{q_{k+1}}$, which are also separated by
at least $2^{-L_{k+1}}$. Let denote by
$C_{q_{k+1}}$ (resp. by $D_{q_k}$) the union of 
$2^{-5(k+1) - 3 + L_{k+1}}$ of these intervals
(resp. of the corresponding $2^{-5k-5 + L_k}$
intervals from $C_{q_k}$).
In particular, we have ${\rm mes} C_{q_k}
= 4 {\rm mes} D_{q_k} = 2^5 {\rm mes} C_{q_{k+1}}$.}

\medskip

We deduce from (4.2), (4.3) and Lemma 2 that
$$
{\rm mes} (B_{q_1}) \ge 1 - \sum_{x=q_0}^{q_1} \, {\rm mes} A(x)
\ge 31/32. 
$$
Consequently, $B_{q_1}$ is the union of at least $2^{L_1-1}$
intervals of length $2^{-L_1}$.
By (4.6), the set $B_{q_0}$ is the union of at least  
$2^{L_0-1}$ intervals of length $2^{-L_0}$.  
This allows us to define the sets $C_{q_0}, D_{q_0}$  
and $C_{q_1}$.  
This proves $(P_0)$. 

Let $k$ be a non-negative
integer such that $(P_k)$ holds, and consider
the set $B'_{q_{k+2}} := C_{q_{k+1}} \cap B_{q_{k+2}}$.
Observe that
$$
B'_{q_{k+2}} = C_{q_{k+1}} \setminus
 \biggl(\, \bigcup_{x=q_{k+1}+1}^{q_{k+2}}  A (x) \biggr),
$$
hence
$$
{\rm mes}  B'_{q_{k+2}} \ge {\rm mes}  C_{q_{k+1}}  - 
\sum_{x=q_{k+1}+1}^{q_{k+2}} {\rm mes} 
\bigl(C_{q_{k+1}}\cap A (x) \bigr).
\eqno (4.8)
$$
By construction,
the set $C_{q_k}$
can be written as a union, say
$$
C_{q_k} = \bigcup_{\nu = 1}^{T_{q_k} } J_\nu,
$$
of $T_{q_k}$ dyadic intervals $J_\nu$ of the form 
$$
\left[ {a \over 2^{L_k}}, {a+1 \over 2^{L_k}}\right], 
\quad a\in \Z_{\ge 0},
$$
where $L_k$ is given by (4.7). Let $x \ge q_k^3$ be an integer. 
Since, by (4.6),
$$
2^{L_k} \le {q_k^2 \log_2^2 q_k \over 4 \eps} \le {q_k^3 \over 2}
\le {x \over 2},
$$
each interval $J_\nu$ 
contains at least the rationals $y/x, (y+1)/x$
for some integer $y$, and we infer from (4.3) that
$$
{\rm mes} (J_\nu \cap A(x)) \le    
 {2^4\varepsilon \over ||\alpha x||\, x\log_2^2 x} 
\times {\rm mes} J_\nu. \eqno (4.9)
$$
Summing (4.9) from $\nu = 1$ to $\nu = T_{q_k}$, we get
$$
 {\rm mes }  ( C_{q_k}  \cap A  (x) ) 
\le {2^4\varepsilon \over ||\alpha x||\, x\log_2^2 x} 
\times {\rm mes} C_{q_k}.  \eqno (4.10)
$$
It then follows from (4.10) that
$$
\eqalign{
{\rm mes}(C_{q_{k+1}}\cap A (x)) & \le {\rm mes}(C_{q_{k}}\cap A (x)) 
\cr & \le {2^4\varepsilon \over ||\alpha x||\, x\log_2^2 x} 
\times {\rm mes} C_{q_{k}} \le
{2^{9} \varepsilon \over ||\alpha x||\, x\log_2^2 x} 
\times {\rm mes} C_{q_{k+1}}. \cr } 
$$
Combined with (4.8) and Lemma 2, this gives
$$
{\rm mes}  B'_{q_{k+2}} \ge ( {\rm mes}   C_{q_{k+1}}) \,  
\biggl( 1- 
\sum_{x=q_{k+1}+1}^{q_{k+2}} \, 
{2^{9} \varepsilon \over ||\alpha x||\, x\log_2^2 x} \, \biggr) \ge
{ {\rm mes}   C_{q_{k+1}} \over 2}.
$$
Thus, at least one quarter of the intervals composing $C_{q_{k+1}}$ 
contains at least $2^{L_{k+2}-L_{k+1} - 2}$
intervals composing $B'_{q_{k+2}}$, 
thus at least $2^{L_{k+2}-L_{k+1} - 3}$
intervals composing $B'_{q_{k+2}}$, if we impose that
these intervals are mutually distant by at least $2^{-L_{k+2}}$.
This allows us to define the sets $C_{q_{k+2}}$ and  
$D_{q_{k+1}}$ with the required properties. 
This proves $(P_{k+1})$.

It then follows that the set
$$
{\cal K} := \bigcap_{k \ge 0} \, D_{q_k} 
$$
is non-empty. 
By construction, every point $\beta$ in this set avoids all the
intervals $E(x, y)$ with $x \ge q_0$, thus, the pair $(\alpha, \beta)$
satisfies (2.1).

To establish that the set ${\cal K}$ has full Hausdorff
dimension, we apply Lemma 1 with
$$
m_k = 2^{L_{k+1}-L_k - 5} \quad
\hbox{and} \quad \eps_k := 2^{- L_{k+1}}.
$$
Note that
$$
{\log (m_1 \ldots m_{k-1}) \over - \log (m_k \eps_k)} 
\ge {\log(32^{-k} 2^{L_k}) \over  - \log (2^{-L_k + 5} )}
$$
We infer from (4.1), (4.5) and (4.7) that
$$
2^{L_k} \ge \delta q_0^{3^k}.
$$
Consequently, 
$$
\lim_{k \to + \infty} \, 
{\log (m_1 \ldots m_{k-1}) \over - \log (m_k \eps_k)} = 1,
$$
and it follows from Lemma 1 that
the set ${\cal K}$ has full Hausdorff dimension.
This completes the proof of our theorem. \cqfd

\vskip 6mm

\centerline{\bf 5. Proofs of Theorems 3, 4, and 5}

\vskip 4mm

The proofs of Theorems 3 and 5 follow exactly the same steps 
as that of Theorem 1. Instead of the intervals
$$
E (x,y) = \biggl[ {y \over x} -{\varepsilon \over ||\alpha x || x^{2}\log_2^2 x}, 
{y \over x} 
+{\varepsilon \over ||\alpha x || x^{2}\log_2^2 x}
\biggr],
$$
we use respectively the intervals
$$
\biggl[ {y \over x} -{\varepsilon \over ||\alpha x || x^{2}
(\log_2 x)^{2-a} \, (\log 1/||\alpha x ||)^a}, 
{y \over x} 
+{\varepsilon \over ||\alpha x || x^{2} 
(\log_2 x)^{2-a} \, (\log 1/||\alpha x ||)^a}
\biggr]
$$
and
$$
\biggl[ {y \over x} -{\varepsilon \over |x|_p x^{2}
(\log_2 x)^{2-a} \, (\log 2/|x|_p)^a}, 
{y \over x} 
+{\varepsilon \over |x|_p x^{2}
(\log_2 x)^{2-a} \, (\log 2/|x|_p)^a}
\biggr].
$$
Furthermore, we apply Lemmas 3 and 4 in place of Lemma 2.

\bigskip

For the proof of Theorem 4, we work directly in the plane.
The idea is the following.
For a triple $(x, y, z)$ of integers and a positive $\eps$, 
the inequality $|x X + y Y + z| \le \eps$ defines
a strip composed of points $(X, Y)$ close to the line
$x X + y Y + z =0$.
Since we are working in the unit square, to a given pair $(x, y)$
of integers corresponds a unique $z$, and the length of the intersection
of the line with the unit square is at most equal to $\sqrt{2}$.
Setting
$$
\eps_{x, y} = {\eps \over |xy| \log^2 |xy|},
$$
for a given (very small) positive $\eps$, the strips 
$|x X + y Y + z| \le \eps_{x, y}$ play the same role as the 
intervals (4.3) in the proof of Theorem 1.

Since, for every large integer $q$, we have
$$
\eqalign{
\sum_{q \le xy \le q^3} \, \eps_{x,y}
& \ll \sum_{x=1}^{q^3} \, 
\sum_{y = \lfloor q/x \rfloor}^{\lfloor q^3 x \rfloor} \,
{\eps \over |xy| \log^2 q} \cr
& \ll \sum_{x=1}^{q^3} \, {\eps \over \log q} \ll \eps, \cr}
$$
the Peres--Schlag method can be applied as in the
proof of Theorem 1. We omit the details.

\vskip 6mm

\centerline{\bf 6. Further results}

\vskip 4mm

We gather in the present section several results that
can be obtained with the same method as in the proof
of Theorem 1.

\bigskip

* A result on lacunary sequences.

\proclaim Theorem 6.
Let $M$ be a positive real number
and $(t_j)_{j \ge 1}$ be a sequence
such that $t_{j+1}/t_j > 1 + 1/M$ for $j \ge 1$.
Let $c$ be a real number with $0 < c < 1/10$. 
Let $\eps$ be a positive real number.
Then, the Hausdorff dimension of the set
$$
\{ \xi \in [0,1] :  \forall n \ge 1, ||\xi t_n || \ge c / (M \log M)\}.
$$
is at least $1-\eps$ if $M$ is sufficiently large.

Theorem 6 complements the results from \cite{PeSc,Mo1}.

\bigskip

* The use of the mass distribution principle
enables us to improve Theorem 1 of \cite{Mo3}.

\proclaim Theorem 7.
Let $C_1, C_2$ and $\gamma$ be positive real numbers.
Let $(t_n)_{n \ge 1}$ be a sequence of real numbers such that
$$
C_1 n^{\gamma} \le t_n \le C_2 n^{\gamma},
\quad \hbox{for $n \ge 1$}.
$$
Then, there exist a positive $C$ and an integer $n_0$ such that
the set
$$
\bigcap_{n \ge n_0} \, \biggl\{ \xi \in \R : 
|| \xi t_n || > {C \over n \log n} \biggr\}  \eqno (6.1)
$$
has full Hausdorff dimension.

It is established in \cite{Mo3} that the 
Hausdorff dimension of the set (6.1) 
is at least $\gamma / (\gamma + 1)$.

As an immediate application, we get that
the set of real numbers $\xi$ for which 
$$
\liminf_{n \to + \infty} \, n (\log n) ||\xi n^2|| > 0
$$
has full Hausdorff dimension.

\bigskip

* We have stated homogeneous statements, but the
method as well allows us to deal with
inhomogeneous approximation.

\bigskip

* By means of dyadic arguments as it was done 
in the preprint \cite{Mo5}, it is possible to 
generalize Lemmas 2 and 3 and, eventually, to
establish the following statement.

\proclaim Theorem 8.
Consider real paremeters $A>1,\, 0<\varepsilon<1$  and $\delta >0$.  
Let $\alpha$ be a badly approximable real number, such that
$$
\inf_{q\ge 1} \, q \, ||q\alpha||\ge \delta>0.
$$
Consider  real-valued  functions $\psi_j, \,\,j=0,1,2,$, defined over the non-negative real numbers, satisfying tho following conditions:
$$
\hbox{$\psi_0(x)>0$ for $x$ large enough and $\Psi_0$ 
is increasing;}  \leqno (i)
$$
$$
\psi_j(0)=0,\,\,\, j=1,2; \leqno (ii)
$$
$$
\hbox{$\psi_j$ increases in some interval 
of the type $ [0,\xi ],\,\, \xi >0$ 
and $ \max_{0\le x\le \xi}\psi_j (x) \le 1$}; \leqno (iii)
$$
$$
\max_{x\in {\bf N}} \,
\psi_0(x)\psi_2(x^{1-A}) \le \varepsilon.  \leqno (iv)
$$
Define $\psi_2^{-1}(x)$ 
to be the inverse function to $\psi_2(x)$ 
so $\psi_2^{-1}(\psi_2(x)) = x$.
Suppose that  
$$
\sup_{X\in {\bf N}}\,\,\,\,
\sum_{X\le \nu<AX}\,\,
\sum_{1\le \mu \le \nu +2-\log_2\delta}\,\,
2^{\nu-\mu }\,\times \, \psi_2^{-1}\left( 
{\varepsilon \over \psi_0(2^\nu)\psi_1(2^{-\mu})} \right)
\le  {1 \over 2^6}.
$$
Take an arbitrary sequence of reals $(\eta_q)_{q \ge 1}$.
Then there exists a real number $\beta $ such that
$$
\liminf_{q\to+\infty} \,
\psi_0(q)\psi_1(||q \alpha||)\psi_2 (||q \beta +\eta_q||) >\varepsilon .
$$

The proof of Theorem 8 follows directly the arguments from \cite{Mo5}.

For a real number $a$ with $ 0\le a<1$, if we put
$$
\psi_0(x) = x\log^{2-a} x,\,\,
\psi_1(x)= x (\log 1/x)^a,\,\,
\psi_2 (x) =x ,\,\,\,\, \eta_q=0,\, q \ge 1,
$$
then Theorem 8 implies that there exists 
a real number $\beta$ such that
$$
\liminf_{q \to + \infty} \, q \cdot 
(\log q)^{2-a} \cdot (\log 1/\Vert q \alpha \Vert)^a
\cdot \Vert q \alpha \Vert \cdot \Vert q \beta \Vert > 0  
$$
a result which corresponds to Theorem 3 (and to Theorem 1 if $a=0$),
with the exception of the assertion on the Hausdorff dimension.
Unfortunately, we cannot put here $a=1$.


\vskip 7mm

\centerline{\bf References}

\vskip 7mm

\beginthebibliography{999}

\bibitem{AdBu06}
B. Adamczewski and Y. Bugeaud,
{\it On the Littlewood conjecture in simultaneous
Diophantine approximation},
J. London Math. Soc. 73 (2006), 355--366.

\bibitem{BeVe07}
V. V. Beresnevich and S. L. Velani, 
{\it A note on simultaneous Diophantine approximation on planar curves},
Math. Ann.  337  (2007), 769--796. 

\bibitem{BeKlMa}
V. Bernik, D. Kleinbock, and G. A. Margulis, 
{\it Khintchine-type theorems on manifolds: the convergence 
case for standard and multiplicative versions},
Internat. Math. Res. Notices  (2001),  453--486.

\bibitem{BuLiv}
Y. Bugeaud,
Approximation by algebraic numbers,
Cambridge Tracts in Mathematics 160,
Cambridge, 2004.

\bibitem{BDM}
Y. Bugeaud, M. Drmota, and B. de Mathan,
{\it On a mixed Littlewood conjecture in Diophantine approximation},
Acta Arith. 128 (2007), 107--124.

\bibitem{BuHaVe}
Y. Bugeaud, A. Haynes, and S. Velani,
{\it Metric considerations concerning the mixed
Littlewood conjecture}.
In preparation.

\bibitem{CaSw} 
J. W. S. Cassels and H. P. F. Swinnerton-Dyer,
{\it On the product of three homogeneous linear forms and indefinite
ternary quadratic forms}, Philos. Trans. Roy. Soc. London,
Ser. A, 248 (1955), 73--96.

\bibitem{EKL}
M. Einsiedler, A. Katok, and E. Lindenstrauss,
{\it Invariant measures and the set of exceptions to the
Littlewood conjecture},
Ann. of Math. 164 (2006), 513--560.

\bibitem{EiKl07}
M. Einsiedler and D. Kleinbock,
{\it Measure rigidity and $p$-adic Littlewood-type problems},
Compositio Math. 143 (2007), 689--702.

\bibitem{Fal90}
K. Falconer, 
Fractal geometry. Mathematical foundations and applications. 
John Wiley \& Sons, Ltd., Chichester, 1990.

\bibitem{Gal62}
P. Gallagher,
{\it Metric simultaneous Diophantine aproximations},
J. London Math. Soc. 37 (1962), 387 -- 390.

\bibitem{Ja28}
V. Jarn\'\i k, 
{\it Zur metrischen Theorie der diophantischen
Approximationen}, Pr\'ace Mat.-Fiz. 36 (1928/29), 91--106.

\bibitem{KuNi}
L. Kuipers and H. Niederreiter,
Uniform distribution of sequences. 
Pure and Applied Mathematics. Wiley-Interscience [John Wiley \& Sons], 
New York-London-Sydney, 1974.

\bibitem{Lit68}
J. E. Littlewood,
Some problems in real and complex analysis.
D. C. Heath and Co. Raytheon Education Co.,
Lexington, Mass., 1968.

\bibitem{BdM03}
B. de Mathan,
{\it Conjecture de Littlewood et r\'ecurrences lin\'eaires},
J. Th\'eor. Nombres Bordeaux 13 (2003), 249--266.

\bibitem{BdMTe}
B. de Mathan et O. Teuli\'e,
{\it Probl\`emes diophantiens simultan\'es},
Monatsh. Math. 143 (2004), 229--245.

\bibitem{Mo1}
N. G. Moshchevitin,
{\it A version of the proof for  Peres-Schlag's theorem on lacunary sequences}.
Available at arXiv: 0708.2087v2 [math.NT] 15Aug2007.

\bibitem{Mo2}
N. G. Moshchevitin,
{\it Density modulo 1 of sublacunary sequences: application of
 Peres-Schlag's arguments}.
Preprint, available at arXiv:  0709.3419v2 [math.NT] 20Oct2007

\bibitem{Mo3}
N. G. Moshchevitin,
{\it On small fractional parts of polynomials},
J. Number Theory 129 (2009), 349--357.

\bibitem{Mo4}
N. G. Moshchevitin,
{\it Towards BAD conjecture}.
Available at arXiv: 0712.2423v2 12Apr2008.

\bibitem{Mo5}
N. G. Moshchevitin,
{\it Badly approximable numbers related to the Littlewood conjecture}. Preprint, available at arXiv: 0810.0777.

\bibitem{Peck}
L. G. Peck,
{\it Simultaneous rational approximations to algebraic numbers},
Bull. Amer. Math. Soc. 67 (1961), 197--201. 

\bibitem{PeSc}
Yu. Peres and W. Schlag,
{\it Two Erd\H os problems on lacunary sequences: chromatic numbers
and Diophantine approximations}.
Available at: arXiv: 0706.0223v1. 

\bibitem{PoVe}
A. D. Pollington and S. Velani,
{\it On a problem in simultaneous Diophantine approximation:
Littlewood's conjecture},
Acta Math. 185 (2000), 287--306.

\endthebibliography

\goodbreak

\vskip 8mm

Yann Bugeaud \hfill Nikolay Moshchevitin

Universit\'e Louis Pasteur \hfill Moscow State University

Math\'ematiques \hfill Number Theory

7, rue Ren\'e Descartes \hfill Leninskie Gory 1

67084 STRASBOURG Cedex (France) \hfill MOSCOW (Russian federation)

\medskip

{\tt bugeaud@math.u-strasbg.fr} \hfill {\tt moshchevitin@rambler.ru}

\bye